\documentclass[11pt,reqno]{amsart}
\usepackage{tikz}
\usepackage{subfig}
\usepackage{epstopdf}
\usepackage{epsfig}
\usepackage{math}
\graphicspath{{./figures/}}
\usepackage{tikz}
\usetikzlibrary{shapes,arrows}
\tikzstyle{decision} = [diamond, draw, fill=blue!20, 
    text width=4.5em, text badly centered, node distance=3cm, inner sep=0pt]
\tikzstyle{block} = [rectangle, draw, fill=blue!20, 
    text width=5em, text centered, rounded corners, minimum height=4em]
\tikzstyle{line} = [draw, -latex']
\tikzstyle{cloud} = [draw, ellipse,fill=red!20, node distance=3cm,
    minimum height=2em]
\tikzstyle{decision} = [diamond, draw, fill=blue!20, 
    text width=4.5em, text badly centered, node distance=3cm, inner sep=0pt]
\tikzstyle{block} = [rectangle, draw, fill=blue!20, 
    text width=5em, text centered, rounded corners, minimum height=4em]
\tikzstyle{line} = [draw, -latex']
\tikzstyle{cloud} = [draw, ellipse,fill=red!20, node distance=3cm,
    minimum height=2em]

\usetikzlibrary{positioning}
\tikzset{main node/.style={circle,fill=blue!20,draw,minimum size=1cm,inner sep=0pt},  }

\begin{document}
\title[Evolutionary dynamics via Optimal transport]{Population games and Discrete Optimal transport}
\author[Chow]{Shui-Nee Chow}
\author[Li]{Wuchen Li}
\author[Lu]{Jun Lu}
\author[Zhou]{Haomin Zhou}
\thanks{This work is partially supported by NSF 
Awards DMS\textendash{}1419027,  DMS-1620345, 
and ONR Award N000141310408.}
\maketitle
\begin{abstract}
 We propose a new evolutionary dynamics for population games with a discrete strategy set, inspired by the theory of optimal transport and Mean field games. The dynamics can be described as a Fokker-Planck equation on a discrete strategy set. The derived dynamics is the gradient flow of a free energy and the transition density equation of a Markov process. Such process provides models for the behavior of the individual players in population, which is myopic, greedy and irrational. The stability of the dynamics is governed by optimal transport metric, entropy and Fisher information. 
\end{abstract}

\section{Introduction}
Population games are introduced as a framework to model population behaviors and study strategic interactions in populations by extending finite player games \cite{nash1950equilibrium, sigmund1999evolutionary, von2007theory}. It has fundamental impact on game theory related to social networks, evolution of biology species, virus and cancer, etc \cite{social, learning, shah2010dynamics,cancer}. Nash equilibrium (NE) describes a status that no player in population is willing to change his/her strategy unilaterally. To investigate stabilities of NEs, evolutionary game theory \cite{nowak2006evolutionary, san2012,  sigmund1999evolutionary} has been developed in the last several decades. People from various fields (economics, biology, etc) design different dynamics, called mean dynamics or evolutionary dynamics \cite{hofbauer2003, san2009},  under various assumptions (protocols) to describe population behaviors. Important examples include Replicator, Best-response, Logit and Smith dynamics \cite{matsui1992best, shah2010dynamics, smith1984}, just to name a few. A special class of games, named potential games \cite{hofbauer1988theory, monderer1996potential,san2010} are widely considered. Heuristically, potential games describe the situation that all players face the same payoff function, called potential. Thus maximizing each player's own payoff is equivalent to maximizing the potential. In this case, NEs correspond to maximizers of the potential, which gives natural connections between mean dynamics and gradient flows obtained from minimizing the negative potential. An important example is the Replicator dynamics, which is a gradient flow of the negative potential in the probability space (simplex) with a Shahshahani metric \cite{akin1979geometry, RM, Shahshahani}.
 
Recently, a new viewpoint has been brought into the realm of population games based on optimal transport, see Villani's book \cite{am2006,vil2008} and mean field games in the series work of Larsy, Lions \cite{MFG,de2014,lasry2007}. The mean field games have continuous strategy sets and infinite players \cite{bt2012, blanchet2014nash}. Each player is assumed to make decisions according to a stochastic process instead of making a one-shot decision. More specifically, individual players change their pure strategies \textit{locally} and
simultaneously in a continuous fashion according to the direction that maximizes their own payoff functions
most rapidly. Randomness is also introduced in the form of white noise perturbation. The resulting dynamics for individual players forms a mean field type stochastic differential equation, whose probability density function evolves according to the Fokker-Planck equation. Here Mean field serves as a mediator for aggregating individual players' behaviors. For potential games \cite{de2014}, Fokker-Planck equations can also be viewed as gradient flows of free energies in the probability space. Here free energy refers to the negative expected payoff added with a linear entropy term, which models risks that players take. Moreover, the probability space is treated as a Riemannian manifold endowed with optimal transport metric \cite{am2006, vil2003, vil2008}.

The aim of this paper is to propose a mean dynamics 
on discrete strategy set, which possesses the same connections as that of mean field games and optimal transport theory. It should be noted that it is not a straightforward task to transform the theory on games with continuous strategy set directly to discrete settings. This is due to the fact that the discrete strategy set is no longer a length space, a space that one can define length of curves, and morph one curve to another in a continuous fashion.
To proceed, we employ key tools developed in \cite{li-theory, li-finite, li-thesis} (Similar topics are discussed in \cite{chow2012, erbar2012ricci, maas2011gradient}).
More specifically, we introduce an optimal transport metric on the probability space of the strategy set. With such metric, we derive the gradient flow of the discrete free energy as mean dynamics.

In detail, consider a population game with finite discrete strategy set $S=\{1,\cdots, n\}$. Denote the set of population state  
\begin{equation*}
\mathcal{P}(S)=\{(\rho_i)_{i=1}^n\in \mathbb{R}^n~:~ \sum_{i=1}^n\rho_i=1\ ,~\rho_i\geq 0\ ,~i\in S\}\ ,
\end{equation*}
and payoff function $F_i\colon \mathcal{P}(S)\rightarrow \mathbb{R}$, for any $i\in S$. The derived mean dynamics is given by
\begin{equation}\label{a1}
\begin{split}
\frac{ d\rho_i}{dt}&=\sum_{j\in N(i)} \rho_j[ F_i(\rho)- F_j(\rho)+\beta(\log\rho_j-\log\rho_i)]_+\\
&-\sum_{j\in N(i)}\rho_{i}[ F_j(\rho)- F_i(\rho)+\beta(\log\rho_i-\log\rho_j)]_+ \ ,\\
\end{split}
\end{equation}
where $\beta\geq 0$ is the strength of uncertainty, $\rho_i(t)$ is the probability at time $t$ of strategy $i\in S$, $[\cdot]_+=\max\{\cdot,0\}$, and $j\in N(i)$ if $j$ can be achieved by players changing their strategies from $i$. We call \eqref{a1} Fokker-Planck equation of a game.

Dynamics \eqref{a1} can be viewed from numerous
perspectives. First of all, if the game under consideration is a potential
game, i.e. games for which there exists a term called potential
$\mathcal{F}~:~\mathcal{P}(S)\rightarrow \mathbb{R}$ such that  
$
\frac{\partial}{\partial\rho_i}\mathcal{F}(\rho)=F_i(\rho)
$,
then equation \eqref{a1} can be seen as the gradient flow of the free energy defined as
$$
-\mathcal{F}(\rho)+\beta\sum_{i=1}^n\rho_i\log\rho_i
$$ on a Riemannian manifold $(\mathcal{P}(S), \mathcal{W})$. Here $\sum_{i=1}^n\rho_i\log\rho_i$ is the discrete entropy term and $\mathcal{W}$ is an optimal transport metric defined on the simplex. Secondly, equation \eqref{a1} can be regarded as the transition function of a
nonlinear Markov process. Such Markov process models
individual player's decision making process, which is \textit{local, myopic, greedy and
irrational}. Locality refers to the behavior that a player only compares his/her current
strategy with neighboring strategies, instead of the entire strategy
set. Myopicity means that a player makes his/her decision solely based on the current
available information. Greediness reflects the behavior that players always
selects the strategy that improves his/her payoff most rapidly at the current
time. Lastly and most importantly, by introducing white noise through the so
called log-laplacian term in \eqref{a1}, the Markov process models players'
uncertainty in the decision-making process. This uncertainty may be due to player
making mistakes or risk-taking behavior. The risk-taking interpretation allows us to define the noisy payoff $\bar F_i\colon\mathcal{P}(S)\rightarrow \mathbb{R}$ for each strategy $i$, 
\begin{equation}\label{noisy-payoff}
\bar F_i(\rho):=F_{i}(\rho)-\beta\log\rho_{i}\ .
\end{equation}
Intuitively, the monotonicity of the $\log$ term implies that the fewer players
currently select strategy $i$, the more likely a player is willing to take risk
by switching to strategy $i$. If the strength of the noise ($\beta$ term) was
sufficiently large, the equilibrium would deviate relatively far from that
without noise. 

Dynamics \eqref{a1} has many appealing features. For potential games, since the dynamics is a gradient flow, the stationary points of the free energy, named Gibbs measures, are equilibria of \eqref{a1}. Their stability properties can also be studied by leveraging two key notions, namely, {\em relative entropy} and {\em relative Fisher information} \cite{Fisher, vil2008}. Through their relations with optimal transport metric, we show that the relative entropy converges to 0 as $t$ goes to infinity, and the solution converges to the Gibbs measure exponentially fast. For general games, \eqref{a1} is not a gradient flow, which may exhibit complicated limiting behaviors including Hopf bifurcations. And the noise level introduces a natural parameter for such bifurcations. 




The arrangement of this paper is as follows. In section \ref{Game}, we give a brief introduction to population games on discrete sets. In section 3, we derive \eqref{a1} by an optimal transport metric defined on the simplex set, and introduce the Markov process associated with \eqref{a1} from the modeling perspective. In section 4, we study \eqref{a1}'s long time behavior by relative entropy and relative Fisher information. In section \ref{examples}, we discuss the application of our dynamics by working on some well-known population games.

\section{Preliminaries}\label{Game}
In this paper we focus on population games. Consider a game played by countable infinity many players. Each player in the population selects a pure strategy from the discrete
strategy set $S=\{1,\cdots, n\}$. The aggregated state of the population can be described by the population state $\rho=(\rho_i)_{i=1}^n\in \mathcal{P}(S)$, where $\rho_i$ represents the proportion of players choosing pure strategy $i$ and $\mathcal{P}(S)$ is a probability space (simplex):
\begin{equation*}\label{probs}
\mathcal{P}(S)=\{(\rho_i)_{i=1}^n\in\mathbb{R}^n~:~\sum_{i=1}^n \rho_i=1\ , ~0\leq \rho_i\leq 1\ ,~i\in S\}\ .
\end{equation*}
The game assumes that each player's payoff is independent of his/her identity (autonomous game). Thus all players choosing strategy $i$ have the continuous payoff function $F_{i}:
\mathcal{P}(S)\rightarrow \mathbb{R}$.  


A population state $\rho^*\in \mathcal{P}(S)$ is a Nash equilibrium of the population game if
\begin{equation*}
\rho_i^*>0 ~\textrm{implies that} ~F_i(\rho^*)\geq F_j(\rho^*)\ ,\quad\textrm{for all $j\in S$\ .}
\end{equation*}
The following type of population games has particular importance, in which NEs enjoys various prominent properties.

{A population game is named a {\em potential game}, if there exists a differentiable potential function $\mathcal{F}: \mathcal{P}(S)\rightarrow \mathbb{R}$, such that $\frac{\partial}{\partial\rho_i}\mathcal{F}(\rho)=F_i(\rho)$, for all $i\in S$.}
It is a well known fact that the NEs of a potential game are the
stationary points of $\mathcal{F}(\rho)$.

\noindent{\em Example:} Suppose that a unit mass of agents are randomly matched to play symmetric normal-form game with payoff matrix $A\in \mathbb{R}^{n\times n}$.
At population state $\rho$, a player choosing strategy $i$ receives payoff equal
to the expectation of the others, i.e. $F_i(\rho)=\sum_{j\in S}a_{ij}\rho_j$. In particular, if the payoff matrix $A$ is symmetric, then the game becomes a potential game with potential function $\mathcal{F}(\rho)=\frac{1}{2}\rho^TA\rho$, since $\frac{\partial}{\partial \rho_i}\mathcal{F}(\rho)=F_i(\rho)$.

Given a potential game with potential $\mathcal{F}$, define the {\em noisy potential} \begin{equation*}
\mathcal{\bar F}(\rho):=\mathcal{F}(\rho)-\beta \sum_{i=1}^n\rho_i\log\rho_i\ ,\quad \beta>0\ ,
\end{equation*}
which is the summation of potential and Shannon-Boltzmann entropy. 
In information theory, it has been known for a long time that the entropy is a way to model uncertainties \cite{Fisher}. In the context of population games, such uncertainties may refer to
players' irrational behaviors, making mistakes or risk-taking behaviors. 
In optimal transport theory, the negative noisy potential is usually called the {\em free energy} \cite{vil2003, vil2008}. 

The problem of maximizing each player's payoff with uncertainties is equivalent to maximizing the noisy potential (minimizing the free energy)
\begin{equation*}
\min\{-\mathcal{\bar F}(\rho)~:~\rho\in\mathcal{P}(S)\}\ .
\end{equation*}
We call the stationary points $\rho^*$ of the above minimization the
discrete Gibbs measures, i.e.
$\rho^*$ solves the following fixed point problem
\begin{equation}\label{gibbs}
\rho_i^*=\frac{1}{K}e^{\frac{F_i(\rho^*)}{\beta}}\ ,~\textrm{for any $i\in S$\ , where}\quad K=\sum_{j=1}^n e^{\frac{F_j(\rho^*)}{\beta}}\ .
\end{equation}

\section{Evolutionary dynamics via optimal transport}\label{derivation}
In this section, we first introduce an optimal transport metric for population games. Based on such a distance, we propose another approach to evolutionary dynamics by optimal transport theory, see references in Villani's book \cite{vil2003, vil2008}. For potential games, such dynamics can be viewed as gradient flows of free energies. 
\subsection{Optimal transport metric for games}
To introduce the optimal transport metric, we start with the construction of strategy graphs. A strategy graph $G=(S,E)$ is a {\em neighborhood} structure imposed on the strategy set $S=\{1,\cdots,
n\}$. Two vertices $i,j\in S$ are connected in $G$ if players who currently choose strategy $i$ is able to switch to strategy $j$. Denote the neighborhood of $i$ by
\begin{equation*}
N(i)=\{j\in S\mid  (i,j)\in E \}\ .
\end{equation*}
 For many games, every two strategies are connected, making $G$ a complete
graph. In other words, $N(i)=S\setminus \{i\}$, for any $i\in S$.
For example,
the strategy set of Prisoner-Dilemma game is either Cooperation (C) or Defection
(D), i.e. $S=\{C, D\}$. Thus, the strategy graph is
   \begin{center}
\begin{tikzpicture}[->,shorten >=1pt,auto,node distance=3cm,
        thick,main node/.style={circle,fill=blue!20,draw,minimum size=1cm,inner sep=0pt]}]
   \node[main node] (1) {$D$};
    \node[main node] (2) [left of=1]  {$C$};
    \path[-]
    (2) edge node {} (1);
          \node[anchor=south] at ( 0,0.5) {$F_D(\rho)$};
      \node[anchor=south] at ( -3,0.5) {$F_C(\rho)$};     
\end{tikzpicture}
\end{center}


For any given strategy graph $G$, we can introduce an optimal transport metric on the simplex $\mathcal{P}(S)$. Denote the interior of $\mathcal{P}(S)$ by $\mathcal{P}_o(S)$. 
 
Given a function $\Phi\colon S\to \mathbb{R}$, define
$\nabla\Phi\colon S\times S\to \mathbb{R}$ as
\begin{equation*}
\nabla\Phi_{ij}=\begin{cases} \Phi_i-\Phi_j\quad &\textrm{if $(i,j)\in E$;}\\
0\quad &\textrm{otherwise}.
\end{cases}
\end{equation*} 
Let $m\colon S\times S\to \mathbb{R}$ be an anti-symmetric flux function such that $m_{ij} = -m_{ij}$. The divergence of $m$, denoted as $\textrm{div}(m)\colon S\rightarrow \mathbb{R}$, is defined by \begin{equation*}
\textrm{div}(m)_i = -\sum_{j\in N(i)}m_{ij}\ .
\end{equation*}
For the purpose of defining our distance function, we will use a particular flux function
\begin{equation*}
m_{ij}=\rho \nabla\Phi:=g_{ij}(\rho)\nabla\Phi_{ij}\ ,
\end{equation*}
where $g_{ij}(\rho)$ represents the discrete probability (weight) on $\mathrm{edge}$ $(i,j)$, defined by 
\begin{equation*}
g_{ij}(\rho)=
\begin{cases}
\rho_j & \bar F_j(\rho)<\bar F_i(\rho)\ ;\\
\rho_i& \bar F_j(\rho)>\bar F_i(\rho)\ ;\\
\frac{\rho_i+\rho_j}{2}& \bar F_j(\rho)=\bar F_i(\rho)\ ,\\
\end{cases}
\end{equation*}
Here $\bar F_i(\rho)=F_i(\rho)-\beta\log\rho_i$, is defined in \eqref{noisy-payoff}.

We can now define the discrete inner product on $\mathcal{P}_o(S)$ of $\nabla\Phi$ 
\begin{equation*}
(\nabla\Phi,\nabla\Phi )_\rho:=\frac{1}{2}\sum_{(i,j)\in E} (\Phi_i-\Phi_j)^2g_{ij}(\rho)\ ,
\end{equation*}
where $\frac{1}{2}$ is applied because each edge is summed twice, i.e. $(i,j)$, $(j, i)\in E$.

The above definitions provide the following distance on $\mathcal{P}_o(S)$. 
\begin{definition}
Given two discrete probability functions $\rho^0$, $ \rho^1\in\mathcal{P}_o(S)$, the Wasserstein metric $\mathcal{W}$ is defined by:
\begin{equation*}\label{metric}
\mathcal{W}(\rho^0,\rho^1)^2=\inf \{\int_0^1(\nabla\Phi,\nabla\Phi)_\rho dt~:~ \frac{d\rho}{dt}+\mathrm{div}(\rho\nabla\Phi)=0\ ,~\rho(0)=\rho^0,~\rho(1)=\rho^1\}\ .
\end{equation*}
\end{definition}
It is known that $(\mathcal{P}_o(S), \mathcal{W})$ is a finite dimensional Riemannian manifold \cite{chow2012, maas2011gradient}. And the metric $\mathcal{W}$ depends on the graph structure of the strategy set.
\subsection{Evolutionary dynamics}\label{derivation}
We shall derive \eqref{a1} as a gradient flow of the free energy on the Riemannian manifold $(\mathcal{P}_o(S), \mathcal{W})$.
\begin{theorem}\label{th12}
Given a potential game with strategy graph $G=(S, E)$, potential $\mathcal{F}(\rho)\in C^2(\mathbb{R}^n)$ and a constant $\beta\geq 0$. Then the gradient flow of free energy
\begin{equation*}
-\mathcal{F}(\rho)+\beta\sum_{i=1}^n\rho_i\log\rho_i
\end{equation*}
 on the Riemannian manifold $(\mathcal{P}_o(S), \mathcal{W})$ is the Fokker-Planck equation 
\begin{equation*}
\begin{split}
\frac{ d\rho_i}{dt}&=\sum_{j\in N(i)} \rho_j[ F_i(\rho)- F_j(\rho)+\beta(\log\rho_j-\log\rho_i)]_+\\
&-\sum_{j\in N(i)}\rho_{i}[ F_j(\rho)- F_i(\rho)+\beta(\log\rho_i-\log\rho_j)]_+ \ ,\\
\end{split}
\end{equation*}
for any $i\in S$.
In addition, for any initial $\rho^0\in\mathcal{P}_o(S)$, there exists a unique solution $\rho(t): [0,\infty)\rightarrow \mathcal{P}_o(S)$. And the free energy is a Lyapunov function. Moreover, if $\rho^{\infty}=\lim_{t\rightarrow \infty }\rho(t)$ exists, $\rho^{\infty}$ is one of the Gibbs measures satisfying \eqref{gibbs}.
\end{theorem}
\begin{remark}
We note that if $\beta=0$ and $G$ is a complete graph, the derived Fokker-Planck equation is the Smith dynamics \cite{smith1984}. 
\end{remark}
\begin{remark}
 The strategy graph $G$ is different from the one in evolutionary graph games studied in \cite{allen2014games, lieberman2005evolutionary, graph}. They mainly consider a spatial space as the graph while our graph relates to the strategy set.
\end{remark}
The proof of Theorem \ref{th12} is shown in \cite{li-theory, li-thesis}, see details there.

We can further extend \eqref{a1} as mean dynamics to model general population games without potential. Although \eqref{a1} can no longer be viewed as gradient flows of any sort in this case, yet it is a system of well defined ordinary differential equations in $\mathcal{P}(S)$. 
\begin{corollary}
Given a population game with strategy graph $G=(S, E)$ and a constant $\beta\geq 0$.  
Assume payoff function $F:~\mathcal{P}(S)\rightarrow \mathbb{R}^n$ are continuous. For any  initial condition $\rho^0\in \mathcal{P}_o(S)$, the Fokker-Planck equation 
\begin{equation*}
\begin{split}
\frac{ d\rho_i}{dt}&=\sum_{j\in N(i)} \rho_j[ F_i(\rho)- F_j(\rho)+\beta(\log\rho_j-\log\rho_i)]_+\\
&-\sum_{j\in N(i)}\rho_{i}[ F_j(\rho)- F_i(\rho)+\beta(\log\rho_i-\log\rho_j)]_+ \ ,\\
\end{split}
\end{equation*}
is a well defined flow in $\mathcal{P}_o(S)$. 
\end{corollary}  
The proof is similar to that of Theorem \ref{th12} and hence omitted. 

It is worth mentioning that, for potential games, there may exist multiple Gibbs measures as equilibria of \eqref{a1}. For non-potential games, there
exist more complicated phenomena than equilibria, for example, invariant sets. We illustrate this by a modified Rock-Scissors-Paper game in Section
\ref{examples}, for which Hopf bifurcation exists with respect to the  parameter $\beta$.  
\subsection{Markov process}\label{Markov_process}
In this subsection, we look at Fokker-Planck equation \eqref{a1} from the probabilistic viewpoint. More specifically, we present
a Markov process whose transition function is given by \eqref{a1}.   
From the modeling perspective, such a Markov process models individual player's
decision process that is myopic, irrational and locally greedy. The Markov process $X_{\beta}(t)$
is defined as
\begin{equation}\label{Markov}
\begin{split}
&\textrm{Pr}(X_\beta(t+h)=j\mid X_\beta(t)=i)\\
=&\begin{cases}
( \bar F_j(\rho)- \bar F_i(\rho))_+h+o(h)\ , \quad&\textrm{if}~ j\in N(i)\ ;\\
1-\sum_{j\in N(i)}(\bar F_j(\rho)-\bar F_i(\rho))_+h+o(h)\ ,\quad &\textrm{if}~ j=i\ ;\\
0\ ,\quad &\textrm{otherwise}\ ,
\end{cases}
\end{split}
\end{equation}
where $\bar F_i(\rho)=F_i(\rho)-\beta\log\rho_i$ and $\lim_{h\rightarrow 0}\frac{o(h)}{h}=0$. It can be easily seen that the probability evolution equation of $X_{\beta}(t)$ is exactly \eqref{a1}.



Process $X_{\beta}(t)$ characterizes players' decision making process. Intuitively,
players compare their current strategies with local
strategy neighbors. If the neighboring strategy has payoff higher than their current payoffs, they switch strategies with probability proportional to the difference between the two payoffs. In addition, $X_\beta(t)$ represents an individual player's irrational behavior. This irrationality may be due to players' mistake or willingness to
take risk. The uncertainly of strategy $i$ is quantified by term
$\log\rho_{i}$. The monotonicity of this term intuitively implies that {\em the fewer
players currently select strategy $i$, the more likely players are willing to
take risks by switching to strategy $i$.} For this interpretation,
we call $F_{i}(\rho)-\beta\log\rho_{i}$ the noisy payoff of strategy $i$, where $\beta$ is the noise level.

\section{Stability via Entropy and Fisher information}
In this section, we discuss the long time behavior of \eqref{a1} for potential games. We shall study the convergence properties of the
dynamics \eqref{a1}. Our derivation depends on two concepts, which are extensions of discrete relative entropy and relative Fisher information \cite{convergence}. They are used to measure the closeness between two discrete measures $\rho$ and $\rho^\infty$, Gibbs measure defined by \eqref{gibbs}.

The first concept is the discrete relative entropy ($\mathcal{H}$)
\begin{equation*}
\mathcal{H}(\rho|\rho^\infty):=\beta(\mathcal{\bar F}(\rho^\infty)-\mathcal{\bar F}(\rho))\ .
\end{equation*}
The other is the discrete relative Fisher information ($\mathcal{I}$)
\begin{equation*}
\mathcal{I}(\rho|\rho^\infty):=\sum_{(i,j)\in E}[(\log\frac{\rho_i}{e^{F_i(\rho)/\beta}}-\log\frac{\rho_j}{e^{F_j(\rho)/\beta}})_+]^2\rho_i\ .
\end{equation*} 

We remark that in finite player games, where the potential is a linear function (non mean-field type), $\mathcal{H}$ and $\mathcal{I}$ coincide with the classical relative entropy (Kullback–Leibler divergence) and relative Fisher information respectively, see \cite{li-finite, li-thesis}.

We shall show that $\mathcal{H}(\rho(t)|\rho^\infty)$ converges to 0 as $t$ goes to infinity. We will also estimate the speed of convergence and characterize their stability properties. Before that, we state a theorem that connects $\mathcal{H}$ and $\mathcal{I}$ via gradient flow \eqref{a1}.
\begin{theorem}\label{H}
Suppose $\rho(t)$ is the transition probability of $X_\beta(t)$ of a potential game.
Then the relative entropy decreases as a function of $t$. In other words, 
\begin{equation*}
\frac{d}{dt}\mathcal{H}(\rho(t)|\rho^\infty)<0\ .
\end{equation*}
And the dissipation of relative entropy is $\beta$ times relative Fisher information  
\begin{equation}\label{Fisher}
\frac{d}{dt}\mathcal{H}(\rho(t)|\rho^\infty)=-\beta\mathcal{I}(\rho(t)|\rho^\infty)\ .
\end{equation}
\end{theorem}
The proof is based on the fact that $\mathcal{H}$ (the difference between noisy potentials) decreases along the gradient flow with respect to time. Namely, 
\begin{equation}\label{lyapunov}
\begin{split}
\frac{d}{dt}\mathcal{H}(\rho|\rho^\infty)=&-\beta\frac{d}{dt}\mathcal{\bar F}(\rho(t))=\beta(\nabla\bar F, \nabla \bar F)_\rho\\
=&\beta\sum_{ (i,j)\in E}[(\bar F_j(\rho)-\bar F_i(\rho))_+]^2\rho_i\\
=&\beta\sum_{ (i,j)\in E}[(\log\frac{\rho_i}{e^{F_i(\rho)/\beta}}-\log \frac{\rho_j}{e^{F_j(\rho)/\beta}} )_+]^2\rho_i\ .\\
\end{split}
\end{equation}
This shows that the noisy potential grows at the rate equal to the relative Fisher information. In other words, the population as a whole always seeks to improve the average noisy payoff at the rate equal to the expected squared benefits.

Based on Theorem \ref{H}, we show that the dynamics converges to the equilibrium exponentially fast. Here the convergence is in the sense of $\mathcal{H}$ going to zero. Such phenomenon is called entropy dissipation.
\begin{theorem}[Entropy dissipation]\label{main-theorem}
Let $\mathcal{F}\in C^2(\mathcal{P}(S))$ be a concave potential function (not necessary strictly concave) for a given game. Then there exists a constant $C=C(p^0,G)>0$ such that 
\begin{equation}\label{exp}
\mathcal{H}(\rho(t)|\rho^\infty)\leq e^{-Ct}\mathcal{H}(\rho^0|\rho^\infty)\ .
\end{equation}
\end{theorem}
The proof of Theorem \ref{main-theorem} is readily available by noticing the fact that 
$$\mathcal{I}(\rho|\rho^\infty)< C\beta\mathcal{H}(\rho|\rho^\infty)\ ,$$
and an application of Grownwall inequality. See details \cite{li-thesis, li-theory}. In fact, the exponential convergence is naturally expected because \eqref{a1} is the gradient flow on a Riemannian manifold $(\mathcal{P}_o(S), \mathcal{W})$. 

In fact, a more precise characterization on the convergence rate $C$ in \eqref{exp} can be made. This characterization enables us to address the stability issues of Gibbs measures. Define
\begin{equation}\label{def}
\lambda(\rho)=\min_{\Phi}~~-\textrm{div}(\rho\nabla\Phi)^T\cdot\textrm{Hess}\mathcal{\bar F}(\rho)\cdot \textrm{div}(\rho\nabla\Phi) \ ,
\end{equation}
where the infimum is among all $(\Phi_i)_{i=1}^n\in \mathbb{R}^n$, such that $(\nabla\Phi, \nabla\Phi)_\rho=1$ and $\textrm{Hess}$ represents the Hessian operator in $\mathbb{R}^n$.

\begin{theorem}[Stability and asymptotic convergence rate]\label{stability}
For a potential game with potential $\mathcal{F}(\rho)\in C^2$. Denote its Gibbs measure $\rho^\infty$ by \eqref{gibbs}. If $\lambda(\rho^\infty)>0$, then $\rho^{\infty}$ is an asymptotic stable equilibrium for \eqref{a1}. In addition, for any sufficiently small $\epsilon>0$, there exists a time $T>0$, such that when $t>T$, \begin{equation*}
\mathcal{H}(\rho(t)|\rho^\infty)\leq e^{-2(\lambda( \rho^{\infty})-\epsilon)(t-T)}
\mathcal{H}(\rho^0|\rho^\infty)\ .
\end{equation*}
\end{theorem}

For more details, see \cite{li-theory}. The above convergence results, including the quadratic minimization \eqref{def}, shares many similar properties with continuous cases. For example, Ricci curvature lower bound and Yano's formula are well defined on discrete strategy set. See details in \cite{li-theory, li-finite, erbar2012ricci, vil2008}.  
\section{Examples}\label{examples}
In this section, we investigate \eqref{a1} by applying it
to several well-known population games.  
\noindent{\em Example 1: Stag Hunt.}  The point we seek to convey in this example is that the noisy payoff reflects the \textit{rationality} of the population. The symmetric normal-form game with payoff matrix 
\begin{equation*}
A=\begin{pmatrix}
h&h\\
0&s
\end{pmatrix}
\end{equation*}
is known as Stag Hunt game. 
Each player in a random match needs to decide whether to hunt for a hare (h) or
stag (s). Assume $s\geq h$, which means that the payoff of a stag is larger than
a hare. This population game has three Nash equilibria: two pure equilibria
$(0,1)$, $(1,0)$, and one mixed equilibrium $(1-\frac{h}{s}, \frac{h}{s})$.

\begin{figure}[H]
 \subfloat[$\beta=5$]
{\includegraphics[scale=0.23]{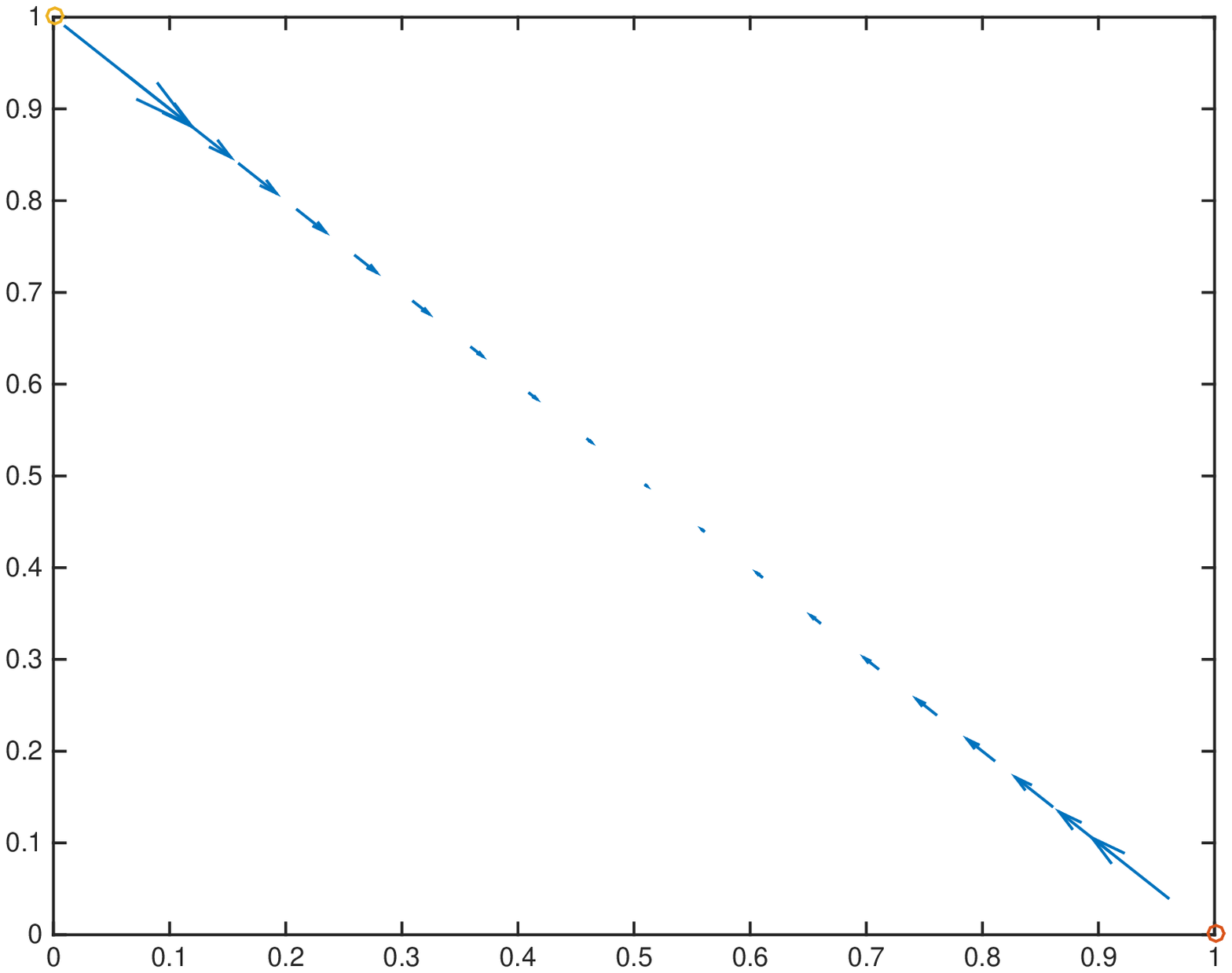}}\hspace{0cm}
 \subfloat[$\beta=0.5$]
{\includegraphics[scale=0.15]{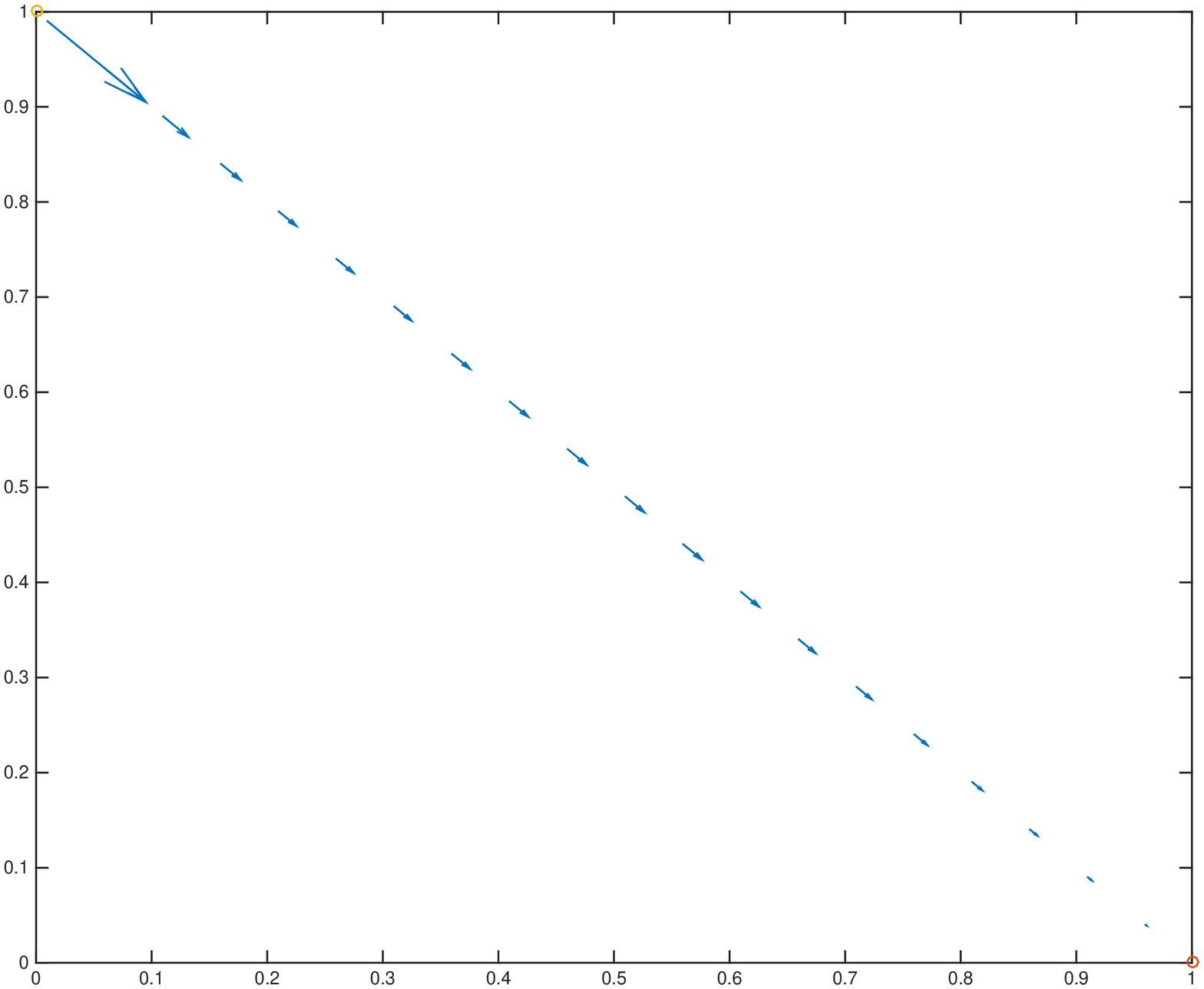}}\\
 \subfloat[$\beta=0.1$]
{\includegraphics[scale=0.23]{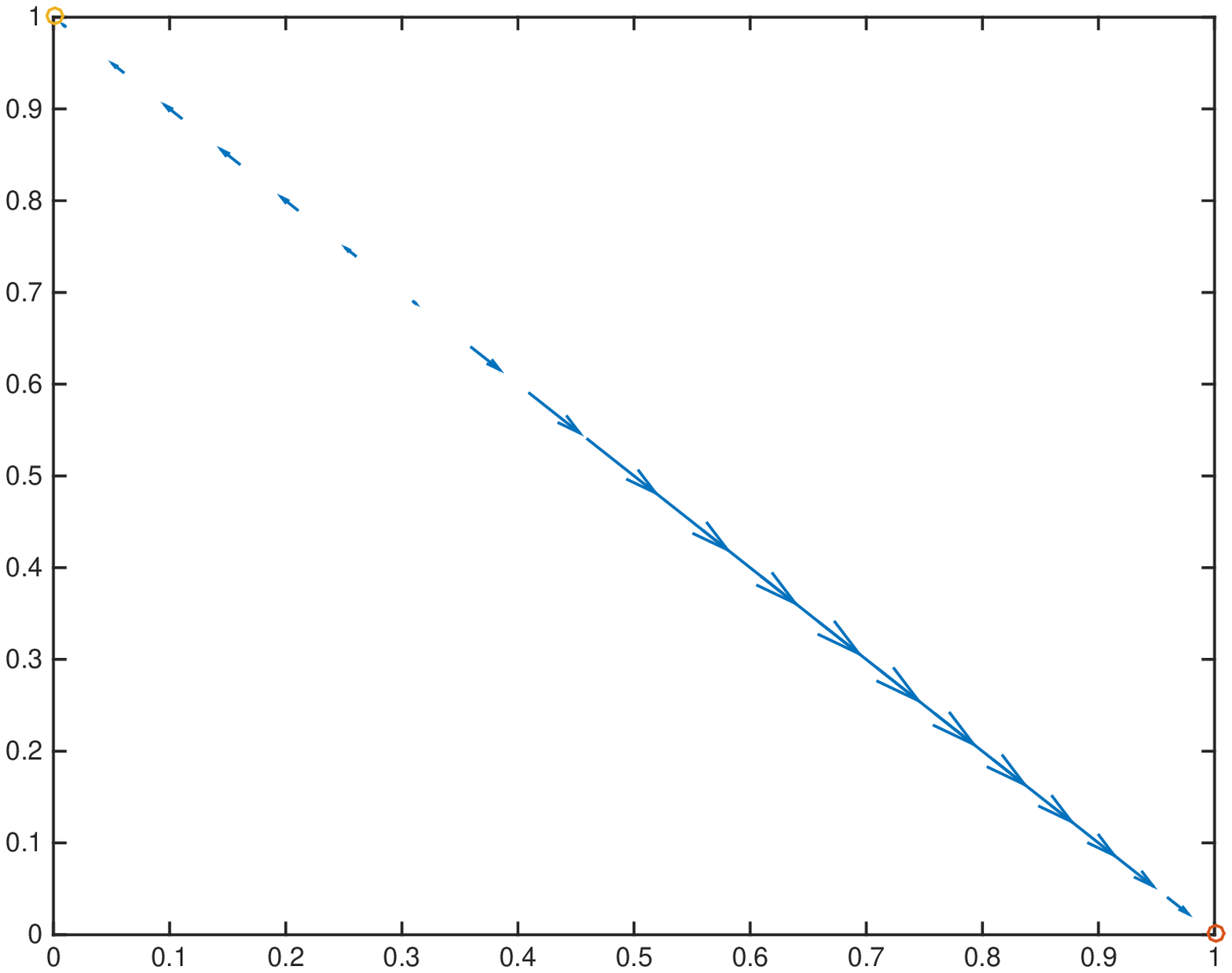}}\hspace{0cm}
 \subfloat[$\beta=0$]
{\includegraphics[scale=0.23]{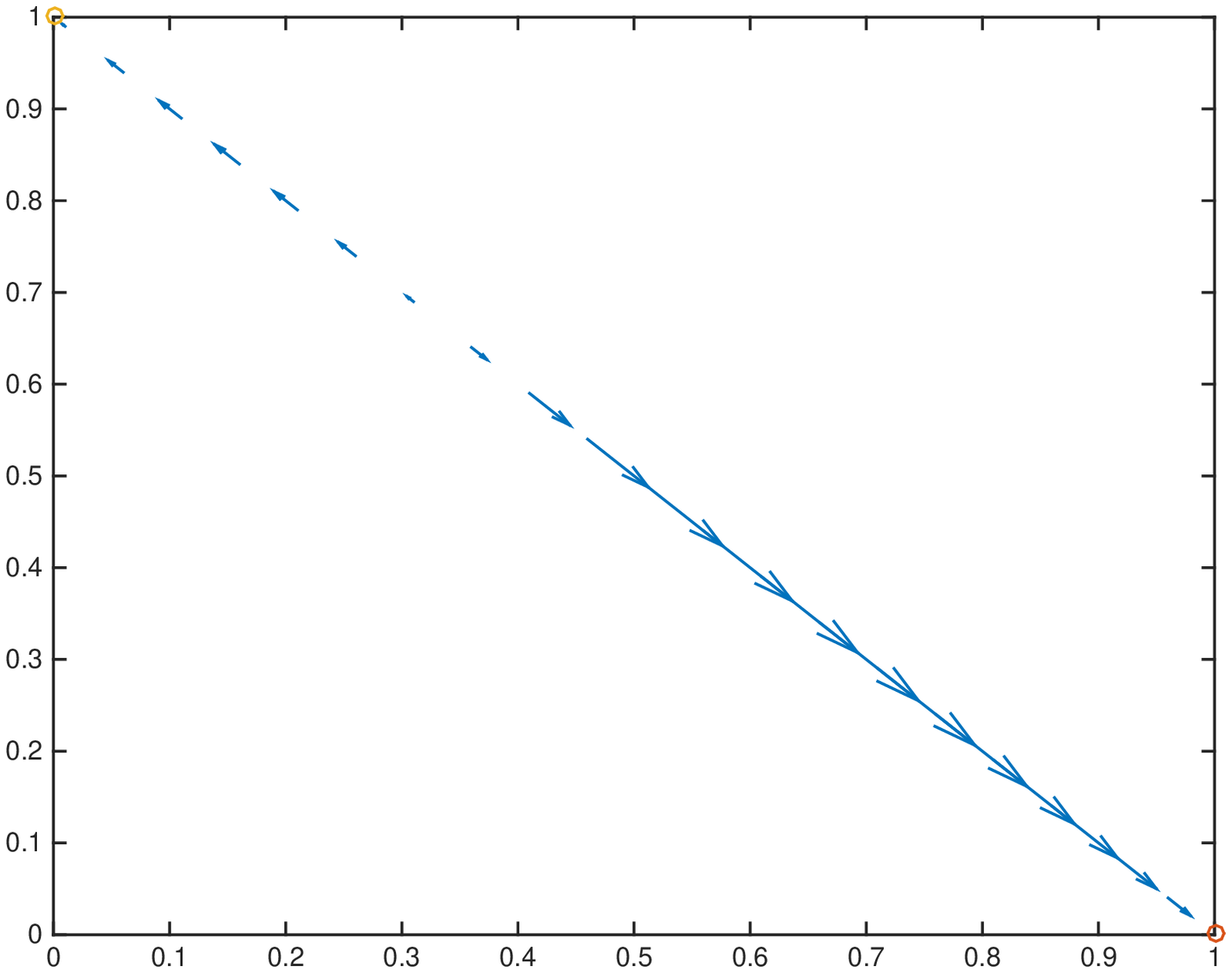}}\\
\caption{Stag and Hare}
\label{stag-hare}
\end{figure}

In particular, let $h=2$ and $s=3$. The population state is $\rho=(\rho_H,\rho_S)^T$ with payoff $F_H(\rho)=2$ and $F_S(\rho)=3\rho_S$.
Then Fokker-Planck equation \eqref{a1} becomes 
\begin{equation*}
\begin{cases}
\dot \rho_H= \rho_S[2-3\rho_S+\beta\log\rho_S-\beta\log\rho_H]_+-\rho_{H}[-2+3\rho_S+\beta\log\rho_H-\beta\log\rho_S]_+\\
\dot\rho_S=\rho_H[3\rho_S-2+\beta\log\rho_H-\beta\log\rho_S]_+-\rho_{S}[-3\rho_S+2+\beta\log\rho_S-\beta\log\rho_H]_+\ .
\end{cases}
\end{equation*}
The numerical results are in Figure \ref{stag-hare}. One can easily see that if
the noise level $\beta$ is 
sufficient small, the perturbation doesn't affect the limit behavior of the mean
dynamics. On the other hand, if noise level $\beta$ is large enough, 
\eqref{a1} settles around $(\frac{1}{2}, \frac{1}{2})$. Lastly, if the noise
level is moderate, Equation \eqref{a1} has  $(1,0)$ as the unique equilibrium.    

The above observation has practical meanings. Namely, if the perturbation is large
enough, it turns out that people always choose to hunt hare (NE $(1,0)$). This
is a safe choice as players can get at least a hare, no matter how the
others behave. This appears even more so  if comparing with the
state $(0,1)$ for which the player receives nothing. If the perturbation is
small and initial population appears to be more cooperative, people will choose
to hunt the stag.  This is a rational move because stag is definitely better
than hare.  

\begin{figure}[H]
 \subfloat[$\beta=0$]
{\includegraphics[scale=0.3]{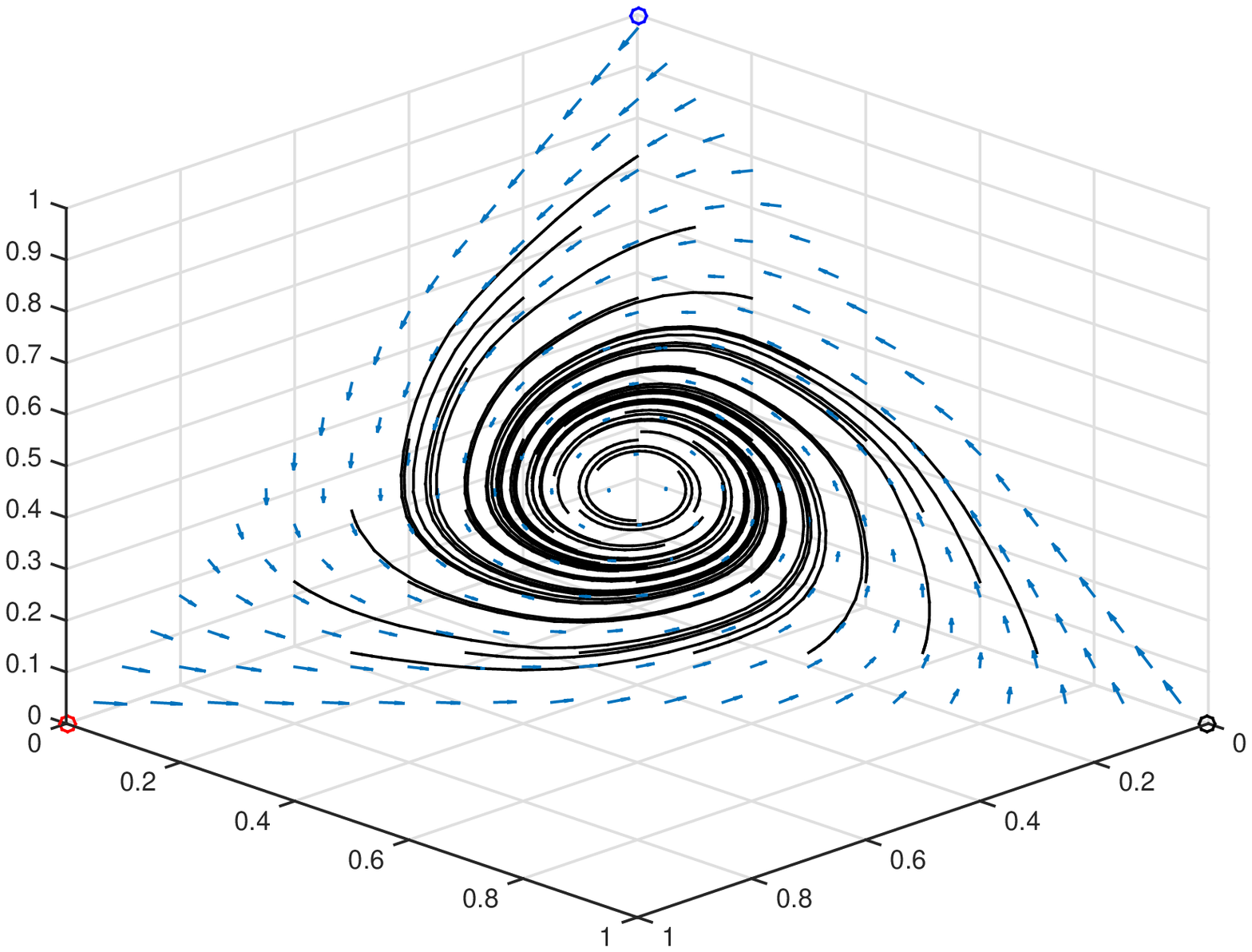}}\hspace{1cm}
 \subfloat[$\beta=0.1$]
{\includegraphics[scale=0.3]{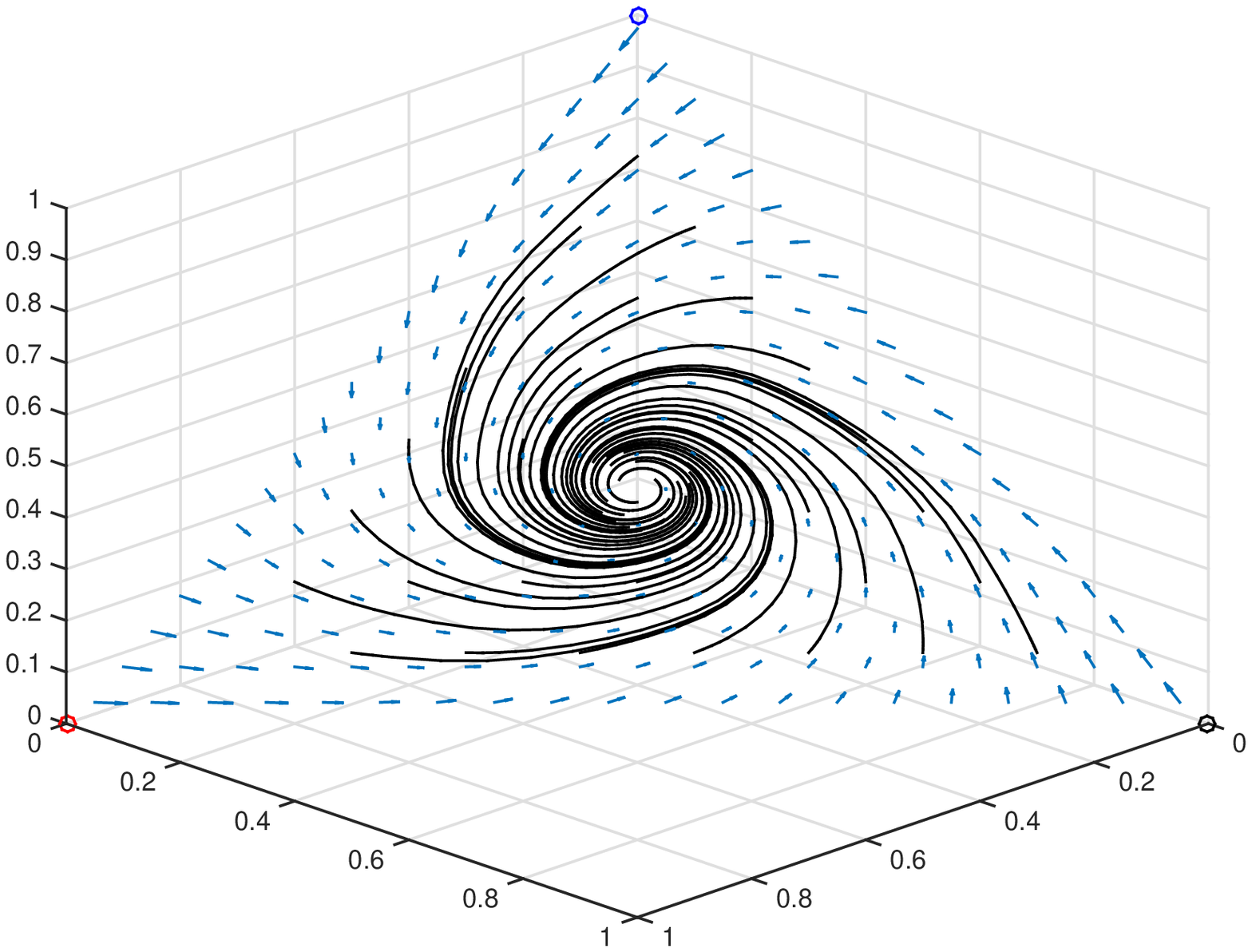}}\\
\caption{Rock-Scissors-Paper}
\label{Rock-Scissors-Paper}
\end{figure}

\noindent{\em Example 2: Rock-Scissors-Paper game}. 
Rock-Scissors-Paper has payoff matrix
\begin{equation*}
A=\begin{pmatrix}
0&1&-1\\
-1&0&1\\
1&-1&0\\
\end{pmatrix}.
\end{equation*}
The strategy set is $S=\{r, s, p\}$. The population state is
$\rho=(\rho_r,\rho_s,\rho_p)^T$ and the payoff functions are $F_r(\rho)=\rho_s-\rho_p$, $F_s(\rho)=-\rho_r+\rho_p$ and $F_p(\rho)=\rho_r-\rho_s$. By solving \eqref{a1}, we find that there is one unique Nash equilibrium  around $\rho^*=(\frac{1}{3},
\frac{1}{3}, \frac{1}{3})$ for various $\beta$s. The result can be found in Figure \ref{Rock-Scissors-Paper}. 

\begin{figure}[H]
 \subfloat[$\beta=0.5$]
{\includegraphics[scale=0.3]{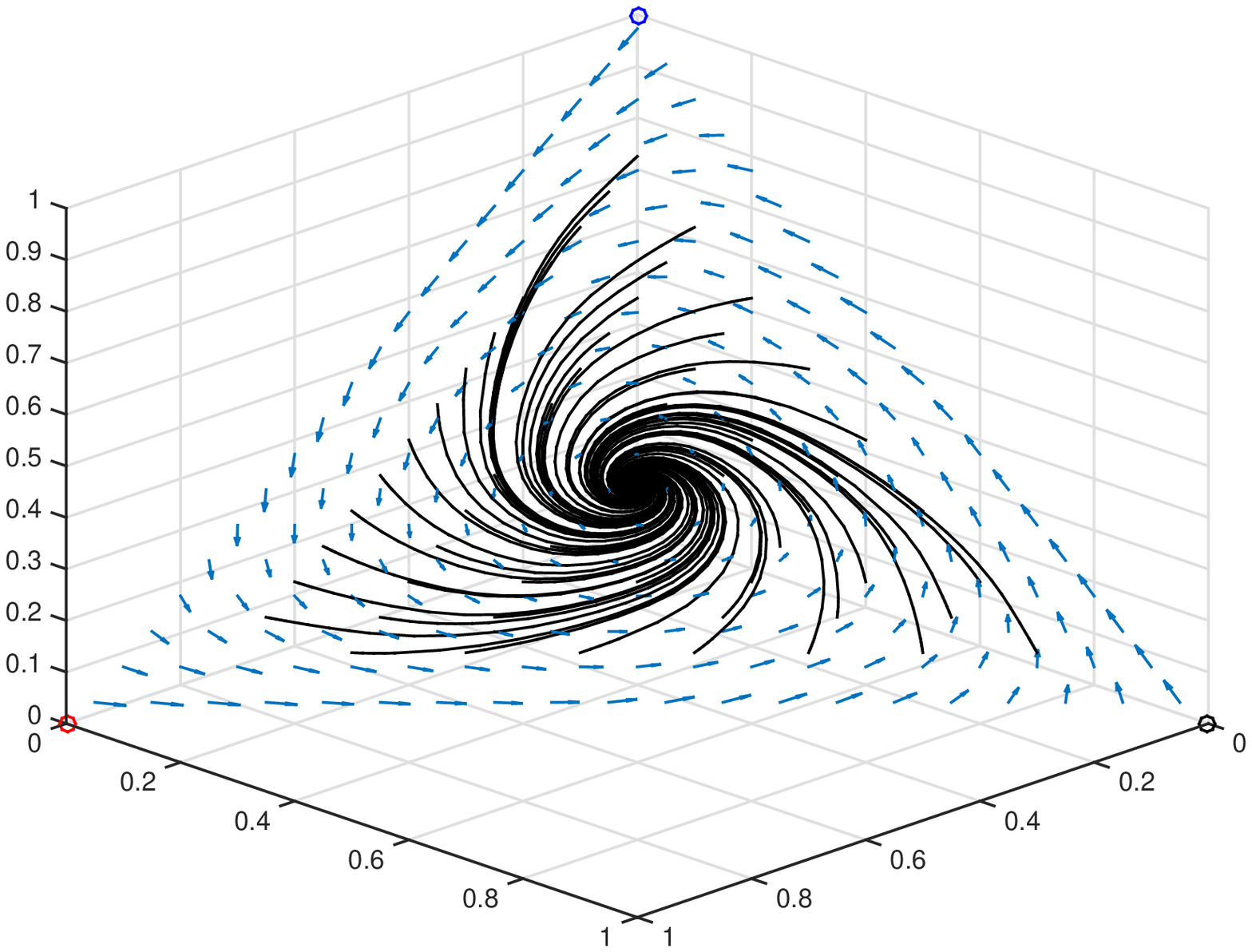}}
 \subfloat[$\beta=0.1$]
{\includegraphics[scale=0.3]{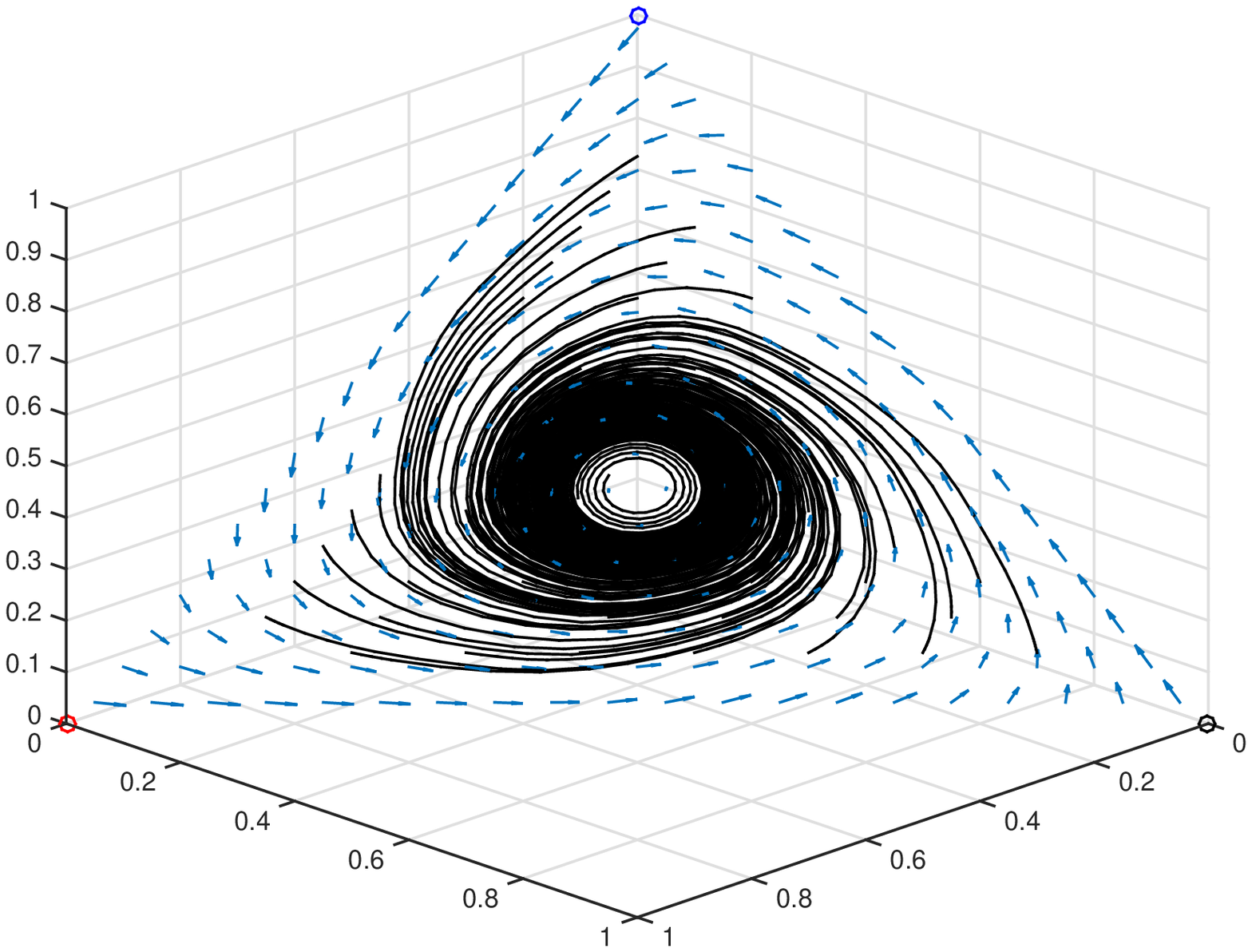}}\\
 \subfloat[$\beta=0$]
{\includegraphics[scale=0.3]{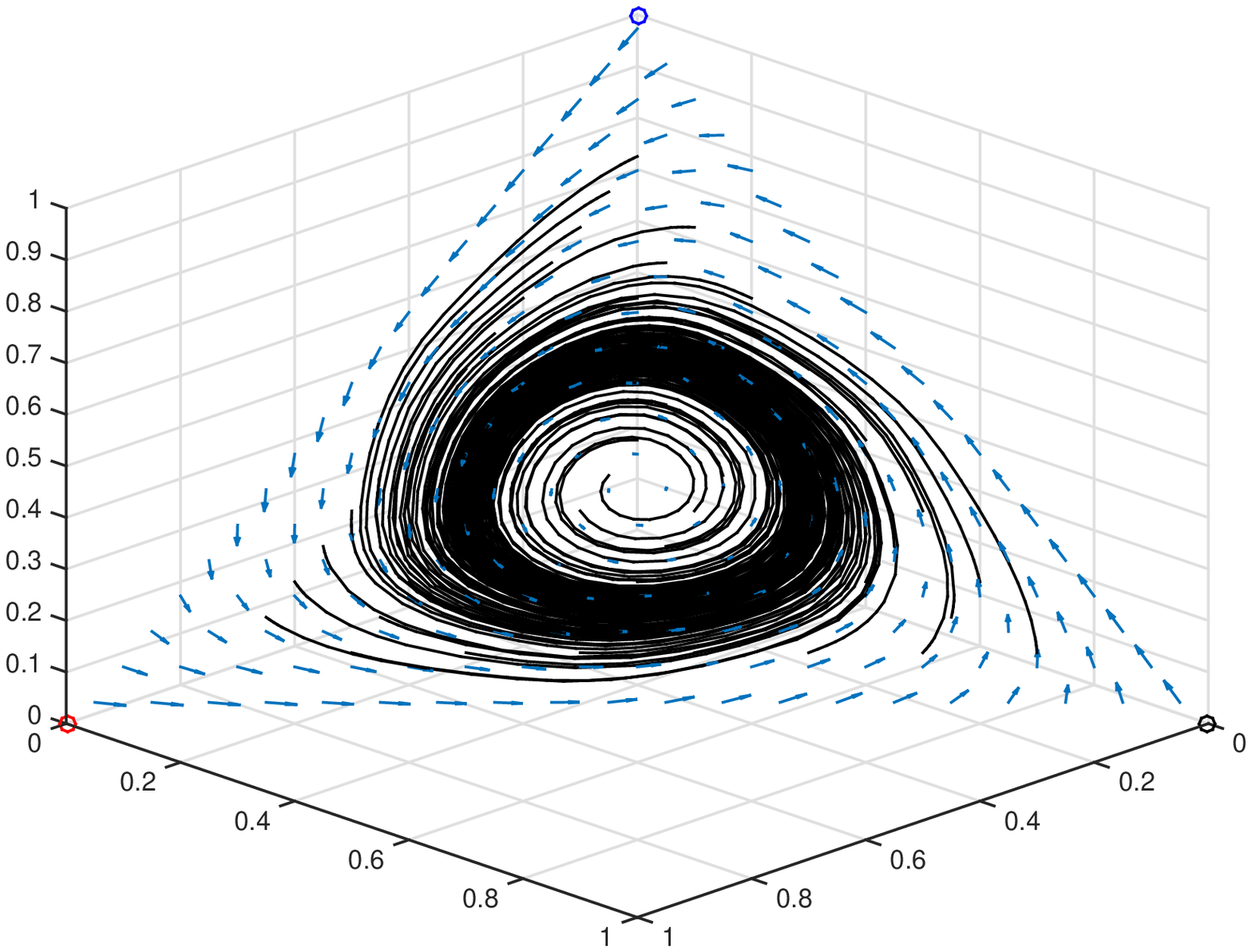}}
\caption{Modified Rock-Scissors-Paper}
\label{Bad-Rock-Scissors-Paper}
\end{figure}

\noindent{\em Example 3.}
We show an example with Hopf Bifurcation. Consider a modified Rock-Scissors-Paper game with payoff matrix
\begin{equation*}
A=\begin{pmatrix}
0&2&-1\\
-1&0&2\\
2&-1&0\\
\end{pmatrix}
\end{equation*}
The strategy set is $S=\{r, s, p\}$. The population state is
$\rho=(\rho_r,\rho_s,\rho_p)^T$ and the payoff functions are
$F_r(\rho)=2\rho_s-\rho_p$, $F_s(\rho)=-\rho_r+2\rho_p$ and
$F_p(\rho)=2\rho_r-\rho_s$. We find that there is Hopf bifurcation for Equation
\eqref{a1}. If $\beta$ is large, there is a unique equilibrium
around $(\frac{1}{3}, \frac{1}{3}, \frac{1}{3})^T$. If $\beta$ goes to $0$,
the solution approaches to a limit cycle. The results are in Figure \ref{Bad-Rock-Scissors-Paper}. 

\begin{figure}[t!]
 \subfloat[$\beta=0$]
{\includegraphics[scale=0.3]{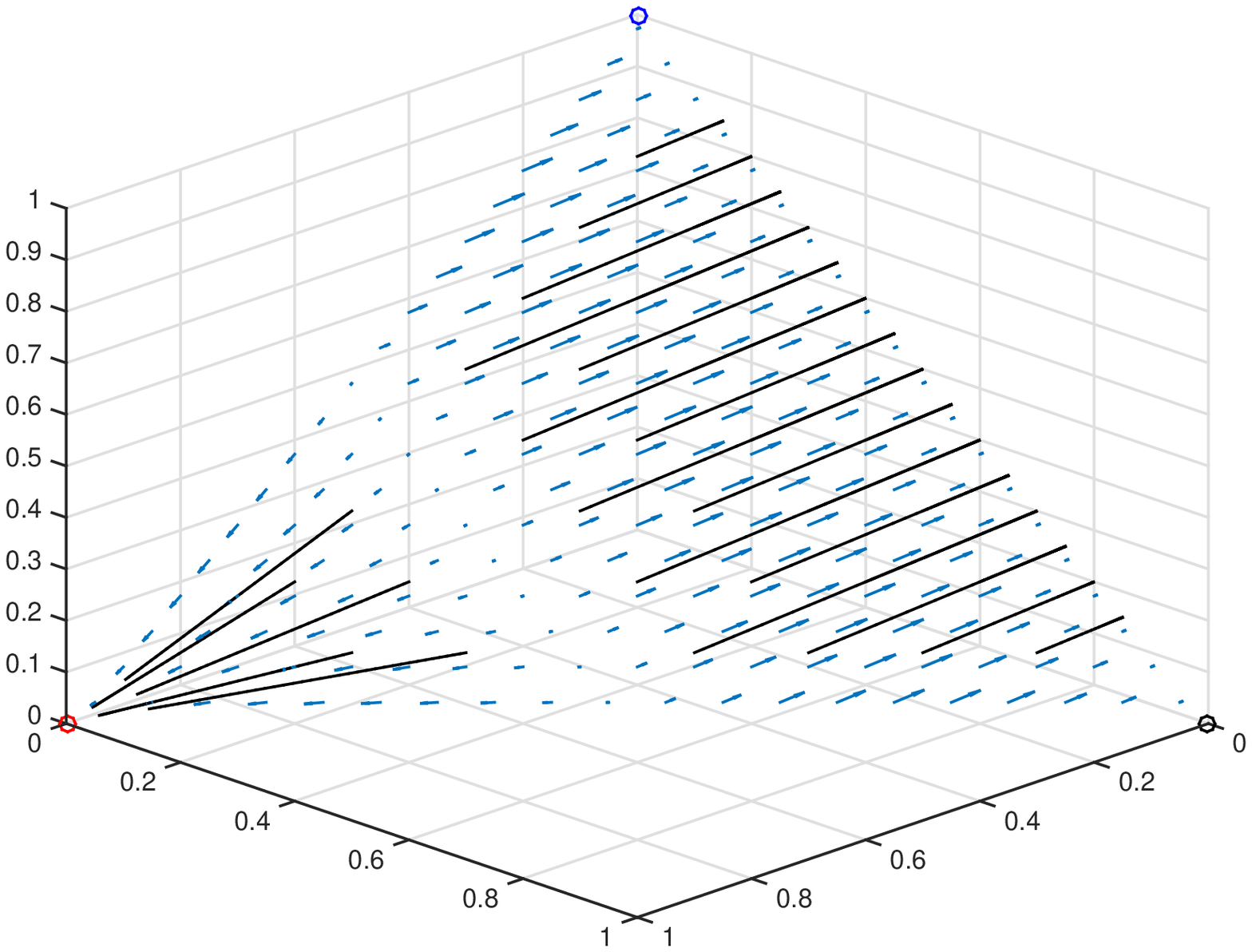}}\hspace{1cm}
 \subfloat[$\beta=0.1$]
{\includegraphics[scale=0.3]{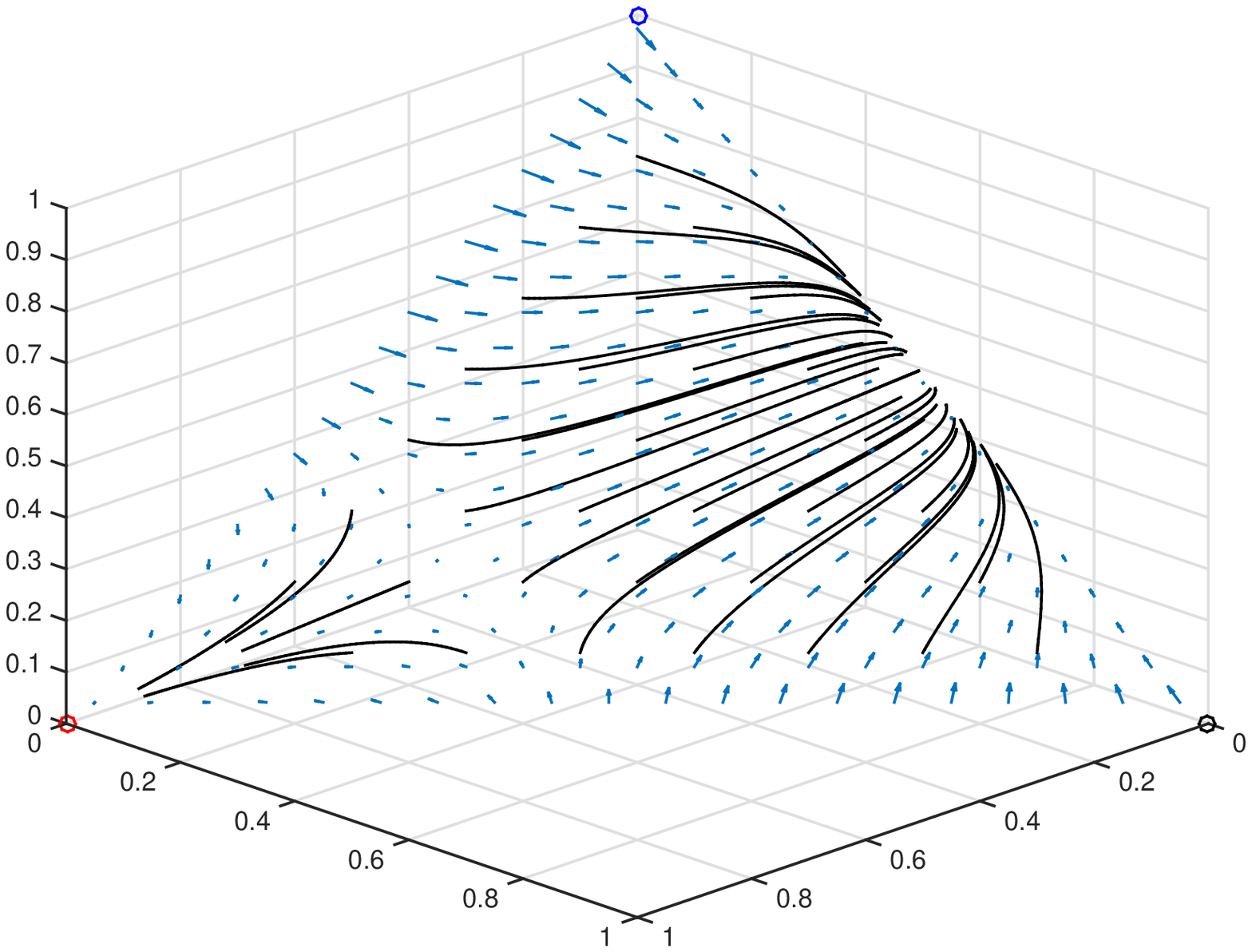}}\\
\caption{Multiple Gibbs measures}
\label{multiple-gibbs}
\end{figure}

\noindent{\em Example 4}. 
We show an example with multiple Gibbs measures. 
Consider a potential game with payoff matrix
\begin{equation*}
A=\begin{pmatrix}
1&0& 0\\
0&1&1\\
0&1&1\\
\end{pmatrix}
\end{equation*}
Denote the strategy set as $S=\{1, 2, 3\}$. The population state is
$\rho=(\rho_1,\rho_2,\rho_3)^T$ and the payoff functions are $F_1(\rho)=\rho_1$,
$F_2(\rho)=\rho_2+\rho_3$ and $F_3(\rho)=\rho_2+\rho_3$. 
We consider three sets of Nash equilibria :   
\begin{equation*}
\{\rho\mid \rho_1=\frac{1}{2}\}\cup \{(1,0,0)\}\cup \{\rho\mid\rho_1=0\}\ ,
\end{equation*}
where the first and third one are lines on the probability simplex $\mathcal{P}(S)$. 
By applying \eqref{a1}, we obtain two Gibbs measures
\begin{equation*}
\{(0,\frac{1}{2},\frac{1}{2})\} \cup \{ (1,0,0)\}
\end{equation*}
as $\beta\rightarrow 0$.
 The vector field is in Figure \ref{multiple-gibbs}.

\begin{figure}[H]
 \subfloat[$\beta=0$]
{\includegraphics[scale=0.3]{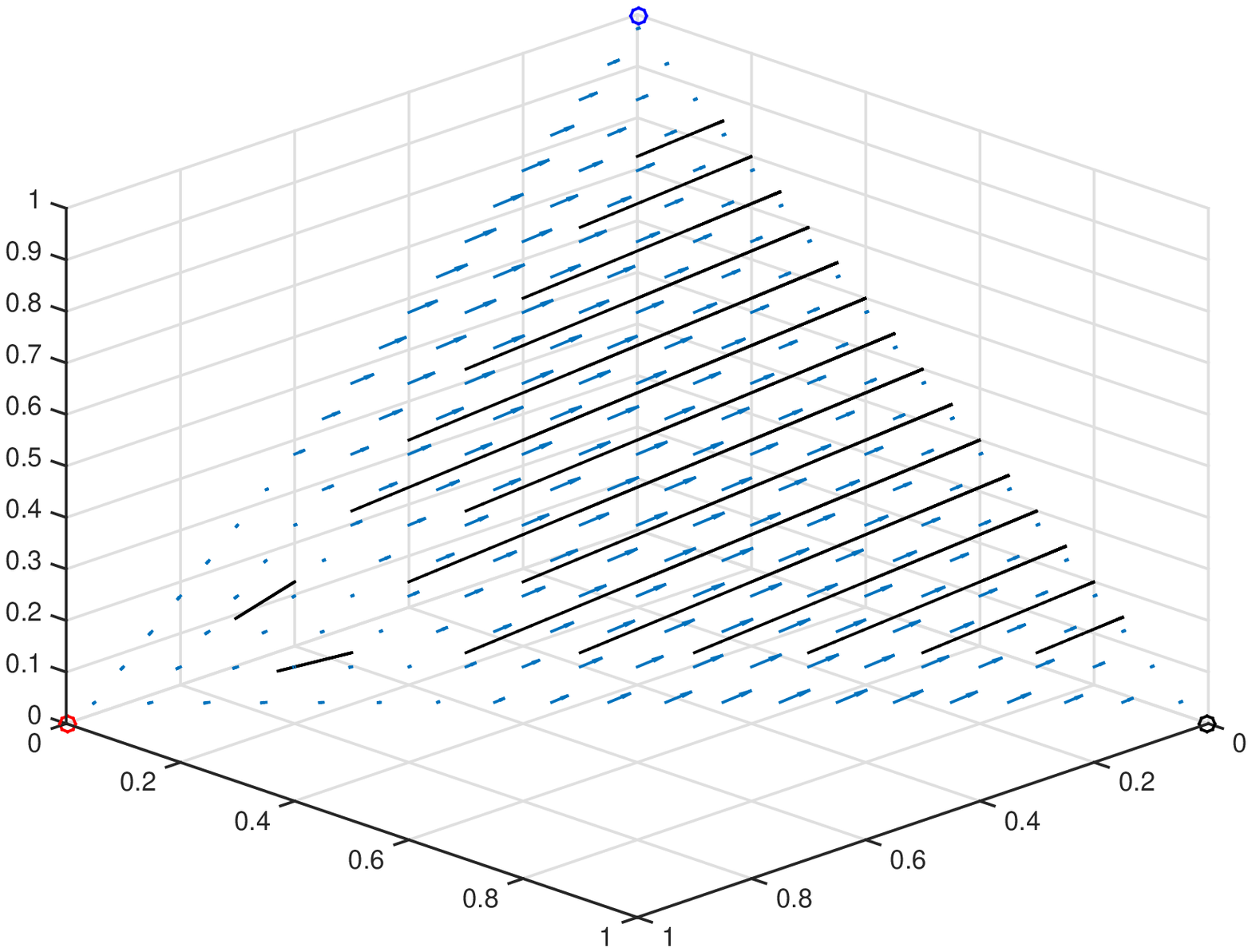}}\hspace{0.45cm}
 \subfloat[$\beta=0.1$]
{\includegraphics[scale=0.3]{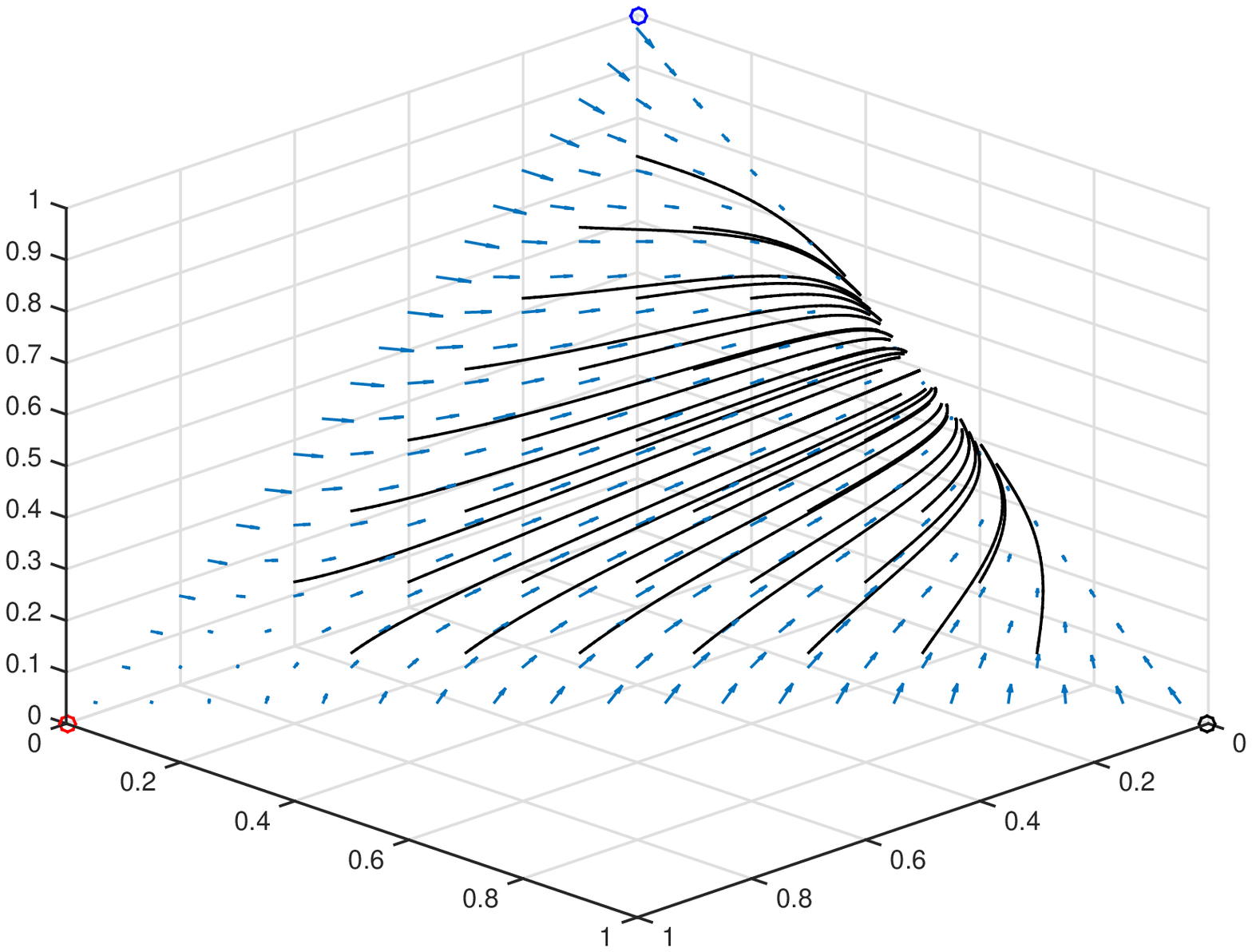}}
\caption{Unique Gibbs measures}
\label{unique-gibbs}
\end{figure}

\noindent{\em Example 5}. 
As a completion, we introduce a game with unique Gibbs measure. 
Let's consider another potential game with payoff matrix
\begin{equation*}
A=\begin{pmatrix}
\frac{1}{2}&0&0\\
0&1&1\\
0&1&1\\
\end{pmatrix}\ .
\end{equation*}
Here the strategy set is $S=\{1, 2, 3\}$, the population state is
$\rho=(\rho_1,\rho_2,\rho_3)^T$ and the payoff functions are $F_1(\rho)=\frac{1}{2}\rho_1$, $F_2(\rho)=\rho_2+\rho_3$ and $F_3(\rho)=\rho_2+\rho_3$.
There are three sets of Nash equilibria 
\begin{equation*}
\{\rho \mid1-\frac{1}{2}\rho_1=\rho_2+\rho_3\}\cup\{(1,0,0)\} \cup \{\rho~|1=\rho_2+\rho_3\}\ ,
\end{equation*}
By applying Fokker-Planck equation \eqref{a1}, we have a unique Gibbs measure
\begin{equation*}
(0,\frac{1}{2},\frac{1}{2})
\end{equation*}
as $\beta\rightarrow 0$. See Figure \ref{unique-gibbs} for the vector fields.

\section{Conclusion}
In this paper, we proposed a dynamics for population games utilizing optimal transport theory and Mean field games. Comparing to existing models, it has the following prominent features.

Firstly, the dynamics is the gradient flow of the noisy potential in the probability space endowed with the optimal transport metric. The dynamics can also be seen as the mean field type Fokker-Planck equations.  

Secondly, the dynamics is the probability evolution equation of a Markov process. Such processes model players' myopicity, greediness and irrationality. In particular, the irrational behaviors or uncertainties are introduced via the notion of noisy payoff. This shares many similarities with the diffusion or white noise perturbation in continuous cases. 

Last but not least, for potential games, Gibbs measures are equilibria of the dynamics. Their stability properties are obtained by the relation of optimal transport metric, entropy and Fisher information. In general, the dynamics may exhibit more complicated limiting behaviors, including Hopf bifurcations.

\textbf{Acknowledgement}: This paper is mainly based on Wuchen Li's thesis. 

\end{document}